\documentclass[12pt]{article}
\usepackage{hyperref}
\usepackage{cmap}                        
\usepackage[cp1251]{inputenc}            
\usepackage[english]{babel}
\usepackage[left=1in,right=1in,top=1in,bottom=1.5in]{geometry} 
\usepackage{amssymb,amsmath, amsthm, amscd,ifthen}
\usepackage{graphicx}
\usepackage[svgnames]{xcolor}
\usepackage{tikz}
\usetikzlibrary{patterns,decorations.text,positioning,arrows,shapes,decorations.markings}
\tikzstyle{vertex}=[circle,draw=black,fill=black,inner sep=0,minimum size=0.2cm,text=white,font=\footnotesize]

\date{}

\title{Tilings with noncongruent triangles}

\author{Andrey Kupavskii\thanks{EPFL, Lausanne and MIPT, Moscow. Supported in part by the grant N 15-01-03530 of the Russian Foundation for Basic Research and by the EPSRC grant no. EP/N019504/1. E-mail: {\tt kupavskii@ya.ru}.} \and J\'anos Pach\thanks{EPFL, Lausanne and R\'enyi Institute, Budapest. Supported by Swiss National Science Foundation Grants 200020-162884 and 200021-165977. E-mail: {\tt pach@cims.nyu.edu}.} \and G\'abor Tardos\thanks{R\'enyi Institute and Central European University, Budapest. Supported by the Cryptography ``Lend\"ulet'' project of the Hungarian Academy
of Sciences and by the National Research, Development and Innovation Office, NKFIH, projects K-116769 and SNN-117879.}}

\date{}

\begin{document}
\maketitle

\begin{abstract}\noindent
We solve a problem of R. Nandakumar by proving that there is no tiling of the plane with pairwise noncongruent triangles of equal area and equal perimeter. We also show that any tiling of a convex polygon with more than three sides with finitely many triangles contains a pair of triangles that share a full side.
\end{abstract}

\section{Introduction}
In his blog, R. Nandakumar~\cite{Na06, Na14} raised several interesting questions on tilings. They have triggered a lot of research in geometry and topology. In particular, he and Ramana Rao~\cite{NaR12} conjectured that for every natural number $n$, any plane convex body can be partitioned into $n$ convex pieces of equal area and perimeter. After some preliminary results~\cite{BBS10, KaHA14} indicating that the problem is closely related to questions in equivariant topology, the conjecture was settled in the affirmative by Blagojevi\'c and Ziegler~\cite{BlZ14} in the special case where $n$ is a prime power. All other cases, including the case $n=6$, are open.
\smallskip

Nandakumar also asked the following interesting variant of the last question. A family of closed triangles that together cover the whole plane is said to form a {\em tiling} if no two of its members share an interior point.
\medskip

\noindent{\bf Problem 1.} (Nandakumar~\cite{Na14}) {\em  Is it possible to tile the plane with pairwise noncongruent triangles of the same area and perimeter?}
\medskip

The aim of this note is to show that the answer to this question is in the negative.
\medskip

\noindent{\bf Theorem 2.} {\em There is no tiling of the plane with pairwise noncongruent triangles of the same area and the same perimeter.}
\medskip

We start with a trivial observation.

\medskip{\bf Observation 3.} {\em If two triangles of the same area and perimeter share a side, then they are congruent.}
\medskip


\noindent Indeed, if two triangles, $xyz$ and $xyz'$, have the same perimeter, then $z$ and $z'$ must lie on a common ellipse with foci $x$ and $y$. On the other hand, if their areas are also the same, the distances of $z$ and $z'$ from the line $xy$ line must be equal. Consequently, if $xyz$ and $xyz'$ do not coincide, one can obtain one from the other by a reflection through the line $xy$, or the midpoint of the segment $xy$, or the orthogonal bisector of this segment.

Therefore, in order to prove Theorem~2, it is sufficient to establish the following result.

\medskip

\noindent{\bf Theorem 4.} {\em Let $\mathcal T$ be tiling of the plane with triangles of unit perimeter, each of which has area at least $\epsilon>0$. Then there are two triangles in $\mathcal T$ that share a side.}
\medskip

In the periodic tiling depicted in Fig.~1, no two triangles share a side. Therefore, Theorem~4 does not hold without the assumption that all triangles have the same perimeter.

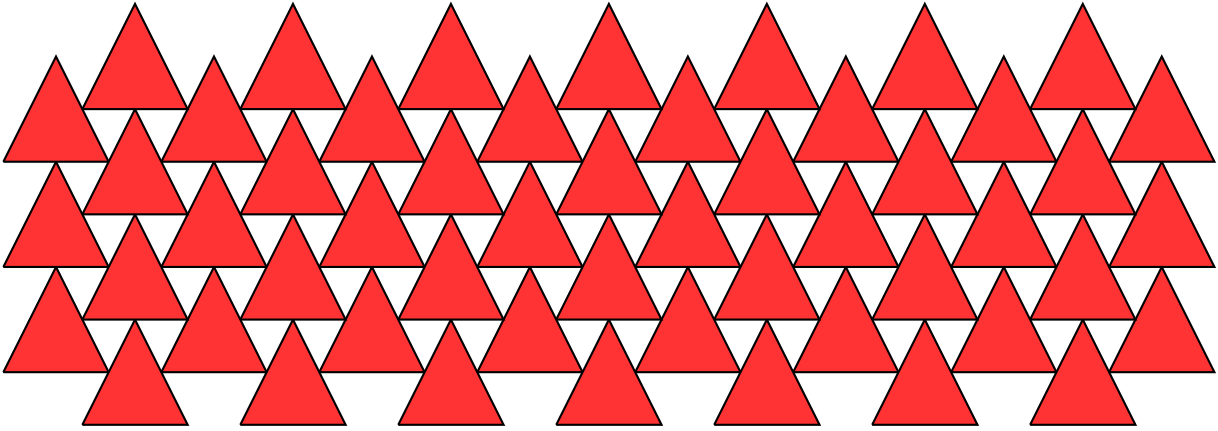
\begin{figure}\centering
\begin{tikzpicture}[scale=1.4]


\foreach \x in {1,...,8}
\foreach \y in {2,...,4}
\filldraw[color=Red!80,thick,draw=black]  ({1.5*\x},{\y}) -- ({1.5*\x+0.5},{\y+1}) -- ({1.5*\x+1},{\y}) -- ({1.5*\x},{\y});
\foreach \x in {1,...,7}
\foreach \y in {1,...,4}
\filldraw[color=Red!80,thick,draw=black]  ({1.5*\x+0.75},{\y+0.5}) -- ({1.5*\x+1.25},{\y+1.5}) -- ({1.5*\x+1.75},{\y+0.5}) -- ({1.5*\x+0.75},{\y+0.5});

\end{tikzpicture}
\caption{Periodic tiling with no two triangles sharing a side.}
\end{figure}

Another way of enforcing that two triangles of a tiling share a side by imposing a condition on the lengths of the sides of the participating triangles. We say that a tiling is {\em locally finite} if any bounded region intersects only a finite number of its members.
\medskip

\noindent{\bf Theorem 5.} {\em Let $\mathcal T$ be a locally finite tiling of the plane with triangles, and suppose that the lengths of their sides belong to interval $[1,2)$. Then there are two triangles in $\mathcal T$ that share a side.}
\medskip

By properly scaling the tiling depicted in Fig.~1, we obtain an example showing that Theorem~5 does not remain true if we replace the interval $[1,2)$ by its closure, $[1,2]$. 





One of the main ideas of the proof of Theorem 2 in a simplified form can be used to obtain a result on triangular tilings of a convex $k$-gon with $k>3$.
\medskip

\noindent{\bf Theorem 6.} {\em Let $k\ge 4$. In any tiling of a convex $k$-gon with finitely many triangles, there are two triangles that share an edge.}
\medskip

We note that, while revising the paper, it was pointed out to us by Roman Karasev and Pavel Kozhevnikov that Theorem~6, for the case $k=4$, was once posed on a Russian mathematics olympiad and was published in Kvant \cite{Kv}.
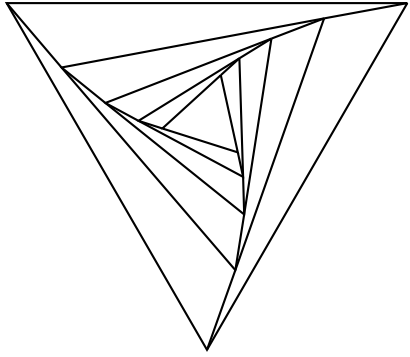
\begin{figure} \centering
\begin{tikzpicture}[scale=0.12]


\foreach \x in {4,...,8}
\draw[thick] ({114.8-10.6*\x}:{(3/2)^\x}) -- ({234.8-10.6*\x}:{(3/2)^\x}) -- ({-5.2-10.6*\x}:{(3/2)^\x}) -- ({114.8-10.6*\x}:{(3/2)^\x});
\foreach \x in {5,...,8}
\draw[thick] ({234.8-10.6*\x}:{(3/2)^\x}) -- ({234.8-10.6*(\x-1)}:{(3/2)^(\x-1)});
\foreach \x in {5,...,8}
\draw[thick] ({114.8-10.6*\x}:{(3/2)^\x}) -- ({114.8-10.6*(\x-1)}:{(3/2)^(\x-1)});
\foreach \x in {5,...,8}
\draw[thick] ({-5.2-10.6*\x}:{(3/2)^\x}) -- ({-5.2-10.6*(\x-1)}:{(3/2)^(\x-1)});

\end{tikzpicture}
\caption{Tiling a triangle with triangles not sharing a side. A triangle is recursively subdivided into $4$ pieces.}
\label{fig3}
\end{figure}

\begin{figure}\centering
\begin{tikzpicture}[scale=0.6]


\draw[thick] (0,0) -- (5,5)-- (10,0)--(5,-5)--cycle;
\draw[thick] (0,0) -- (10,0);
\draw[thick,red] (7.5,0) -- (7.5,-2.5)--(2.5,-2.5)--(2.5,0);
\draw[ultra thick] (6.25,0) -- (7.5,-1.25)--(5,-2.5)--(2.5,-1.25)--(3.75,0);
\draw[thick,red] (5.625,0) -- (6.875,-0.625)--(6.25,-1.875)--(3.75,-1.875)--(3.125,-0.625)--(4.375,0);


\node[fill=white, circle, inner sep = -0pt, minimum size=0pt] at (5,-1.3) {.};
\node[fill=white, circle, inner sep = -0pt, minimum size=0pt] at (5,-1.45) {.};
\node[fill=white, circle, inner sep = -0pt, minimum size=0pt] at (5,-1.6) {.};

\node[fill=black, circle, inner sep = -0pt, minimum size=6pt] at (5,0) {};

\end{tikzpicture}
\caption{Tiling of a square with infinitely many triangles, no pair of which share a side.}\label{figsquare}
\end{figure}
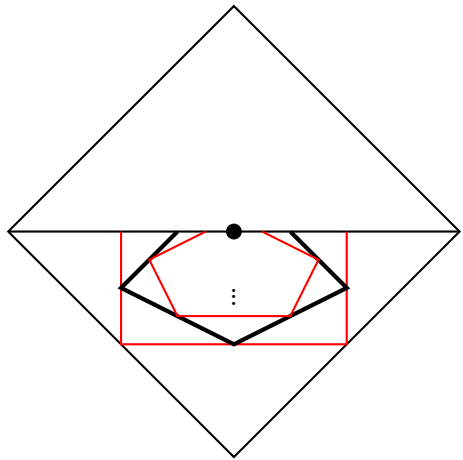

The tiling depicted in Figure~\ref{fig3} shows that Theorem 6 is false for $k=3$. The tiling in Figure~\ref{figsquare} shows that Theorem 6 does not hold for infinite tilings.

The rest of this note is organized as follows.
In Section~2, we establish Theorem~6. The proof of Theorem~4 (which implies Theorem 2) is similar, but more involved. It requires careful asymptotic estimates of certain parameters and is presented in a separate section (Section~3). The proof of Theorem~5 is given in Section~4.
\smallskip

For many other problems on tilings, consult~\cite{BMP05} and~\cite{Na14}.


\section{Tiling a convex polygon---Proof of Theorem 6}

Consider a convex $k$-gon $P$ which is tiled with $t$ triangles such that no two of them share a side. In what follows, by a {\em side} we will always mean a side of one of the $t$ triangles, so that the total number of sides is $3t$. We use the term {\em vertex} for the vertices of the triangles (including the vertices of $P$). Define a simple graph $G$ on this vertex set by connecting two vertices with an edge if the segment between them belongs to the side of a triangle and contains no other vertex.

Obviously, $G$ is a connected  planar graph with $f(G)=t+1$ faces. The number of vertices of $G$ is $v(G)=v_{\rm bd}+v_{\rm int}$, where $v_{\rm bd}$ and $v_{\rm int}$ denote the number of vertices lying on the boundary of $P$ and in the interior of $P$, respectively, and $v_{\rm bd}\ge k$. By Euler's polyhedral formula, the number of edges of $G$ satisfies
$$e(G)=v(G)+f(G)-2=v_{\rm bd}+v_{\rm int}+t-1.$$
The number of edges along the boundary of $P$ is $v_{\rm bd}$, and every edge in the interior of $P$ belongs to the side of precisely two triangles. Denote by $v^*_{\rm int}\le v_{\rm int}$ the number of vertices in the interior of $P$ that subdivide (that is, lie in the interior of) a side of a triangle. In fact, each such vertex subdivides precisely one side. Double counting the edges, we obtain
$$2e(G)=3t+v^*_{\rm int}+v_{\rm bd}.$$
Comparing the last two equations gives
\begin{equation}\label{eq1}
v_{\rm bd}+2v_{\rm int}-v^*_{\rm int}=t+2.
\end{equation}

A {\em minimal} segment that can be decomposed into sides of triangles in two different ways is called a {\em stretch}. The total number of triangles whose sides participate in such a decomposition (that is, lie along the stretch) is called the {\em size} of the stretch. Since no two triangles share a side, the size of every stretch is at least 3. If a stretch has size $s$, then there are precisely $s-2\ge {s\over 3}$ vertices in its interior that {\em subdivide} a side.  Observe that each of the $v^*_{\rm int}$ vertices in the interior of $P$ that subdivide a side lies in the interior of precisely one stretch.

\begin{figure}\centering
\begin{tikzpicture}[scale=3]


\draw[thick] (0,2) -- (4,2) -- (4,0) -- (0,0) -- (0,2);
\draw[thick] (4,2) -- (2.5,0.5);
\draw[thick] (4,0) -- (3,1);
\draw[ultra thick, red] (0,0) -- (2.5,0.5);
\draw[thick] (2.5,0.5) -- (3.34,0.67);
\draw[thick] (0,2) -- (1.5,0.3);
\draw[thick] (4,0) -- (2.5,0.5);
\draw[thick] (4,2) -- (2,0.4) -- (0,2);
\draw[thick] (1,1.2) -- (3.2,1.38) -- (0,2);

\node[fill=black, circle, inner sep = 0pt, minimum size=7pt] at (2.5,0.5) {};
\node[fill=Grey, circle, inner sep = 0pt, minimum size=7pt] at (2,0.4) {};
\node[fill=Grey, circle,inner sep = 0pt, minimum size=7pt] at (1.5,0.3) {};

\end{tikzpicture}
\caption{Tiling of a rectangle with triangles. The thick segment is a stretch of size $4$ with two subdividing vertices marked by grey dots. The only interior vertex which is not a subdividing vertex is marked by a black dot. }
\end{figure}
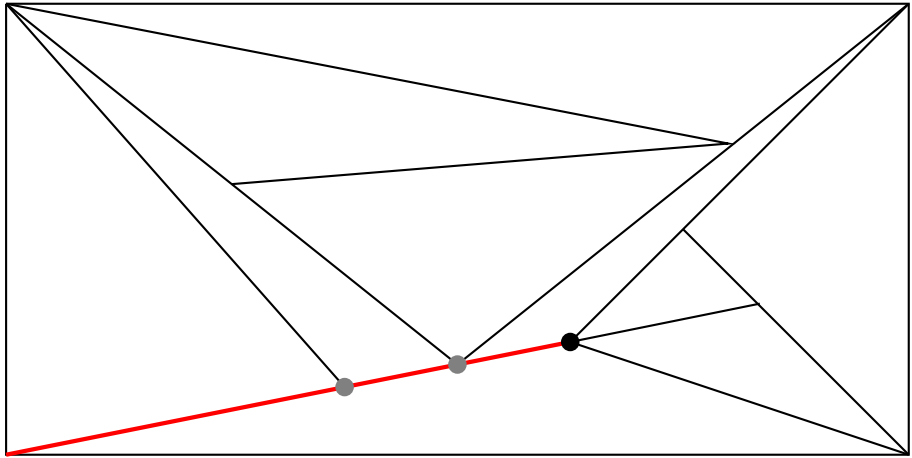

Let $\Sigma$ denote the sum of the sizes of all stretches. Since every side of a triangle, except for the $v_{\rm bd}$ sides that lie on the boundary of $P$, contributes precisely one to this sum (and those on the boundary contribute zero, due to convexity), we have that $\Sigma=3t-v_{\rm bd}$. Along each stretch, the number of subdividing vertices is at least one-third of the size of the stretch. Therefore, the total number of subdividing vertices is at least ${\Sigma\over 3}$. That is,
$$v^*_{\rm int}\ge {3t-v_{\rm bd}\over 3} = t-{v_{\rm bd}\over 3},$$
whence
$$t\le v^*_{\rm int} + {v_{\rm bd}\over 3}.$$
Plugging this inequality into (\ref{eq1}), we obtain that
$$v_{\rm bd}+3(v_{\rm int}-v^*_{\rm int})\le 3.$$
This can hold only if $v_{\rm bd}=3$ and $v_{\rm int}=v^*_{\rm int}$. In particular, $k\le v_{\rm bd}$ must also be equal to $3$. This completes the proof of Theorem~6. $\Box$
\medskip

The above proof also shows that if $P$ is a triangle ($k=3$), then any tiling of $P$ with triangles such that no two triangles share a sides satisfies that

{\bf (i)} there is no vertex subdividing any side of $P$,

{\bf (ii)} every vertex in the interior of $P$ must subdivide a side, and

{\bf (iii)} the size of every stretch is exactly $3$.

\noindent
Despite all these limitations, for any $k>1$, a triangle can be cut into $3k+1$ smaller triangles in several different ways without having two pieces that share a side. For example, this can be achieved by recursively subdividing one of the triangles into $4$ pieces; see, e.g., Fig.~\ref{fig3}.

\section{Equal perimeter tilings---Proof of Theorem 4}

Here we extend the ideas of the proof of Theorem~6 given in the previous section from a finite tiling of a polygon to tilings of the whole plane. We start with a quick overview. Assume for a contradiction that there is a tiling with no pair of triangles sharing a side which satisfies the conditions in Theorem~4. We want to apply the computation described in the previous section to a large but bounded portion of this tiling. In general, we will not be able to apply Theorem~6 directly, because there is no guarantee that there exists a large part of the tiling whose union is a convex polygon. For this reason, as a ``boundary effect'', some small error terms will creep into our calculations. Nevertheless, we will be able to conclude that, analogously to condition (iii) at the end of the last section, {\em most} stretches in the tiling have size $3$. Out of the three sides along a stretch of size 3, one is as long as the other $2$ combined. We call the former side {\em long}, while the other two sides along the same stretch are called {\em short}. Of course, there may be a short side in one stretch that is much longer than a long side in another. Overall, there are twice as many short sides as long sides, but their total lengths are exactly the same. We deduce a strong triangle inequality for the short and long sides of a single triangle. Summing up these inequalities for all triangles, and combining the obtained bound with the above facts, we will arrive at the desired contradiction. Next, we spell out the details of the proof.
\medskip

\noindent{\bf Proof of Theorem~4}. Assume for a contradiction that there is a tiling ${\mathcal T}_{\infty}$ of the plane satisfying the conditions in Theorem~4, but not having two triangles that share a side. Since every triangle has area at least $\epsilon>0$ and unit perimeter, the tiling must be locally finite. The perimeter of each triangle is $1$, which implies that the areas of the triangles are also bounded from above, by the constant $\delta=\sqrt3/36$. Moreover, each side of a triangle has to be longer than $\epsilon'=4\epsilon$, otherwise, its area would be smaller than $\epsilon$. At the end of the proof, we will use the following strong triangle inequality: the total length of any two sides of a triangle in our tiling exceeds the length of the third side by at least some fixed $\epsilon''>0$. The existence of such $\epsilon''$ follows from a compactness argument.

For the rest of the proof, we think of $\epsilon$, $\delta$, $\epsilon'$, and $\epsilon''$ as fixed positive constants. Their exact values are not relevant for our argument.

\smallskip

Let $r>0$ be a sufficiently large number to be specified later, and choose an open circular disk $D(r)$ of radius $r$.  Let ${\mathcal T}_0$ denote the family of all triangles in $T\in {\mathcal T}_{\infty}$ that have a nonempty intersection with $D(r)$.  The union of these triangles is a connected set, but not necessarily {\em simply} connected. Filling up the possible holes with other triangles in ${\mathcal T}_{\infty}$, we turn this set into a simply connected polygonal region $P$. Let ${\mathcal T}$ denote the family of all triangles $T\in {\mathcal T}_{\infty}$ that belong to $P$. The union of these triangles is $P$. Let ${\mathcal T}'\subset{\mathcal T}_{\infty}$ denote the family of triangles that touch the boundary of $P$ from outside.

Since the diameter of every triangle is at most $1/2$, all members of ${\mathcal T}'$ lie in the annulus $D(r+1)\setminus D(r)$, where $D(r+1)$ is the disk of radius $r+1$ concentric with $D(r)$. Using the above lower and upper bounds on the areas of the triangles, and denoting the number of triangles in $\mathcal T$ (and ${\mathcal T}'$) by $t$ (and $t'$, respectively), we have
\begin{equation}\label{1}
\begin{split}
t  &= |\mathcal T| \ge \frac{r^2\pi}\delta=\Omega(r^2),\\
t' &= |{\mathcal T}'| \le \frac{(r+1)^2\pi-r^2\pi}\epsilon=O(r).
\end{split}
\end{equation}

Here we use the asymptotic notation with respect the choice of $r$.
\smallskip

When speaking of {\em vertices} and {\em sides}, we mean the vertices and sides of the triangles in $\mathcal T$. As no two triangles share a side, the total number of sides is $3t$. Define a simple graph $G$ on the set of vertices by connecting two vertices by an {\em edge} if the segment between them belongs to a side and contains no other vertex.

Obviously, $G$ is a connected  planar graph with $f(G)=t+1$ faces, including the exterior face. Let $v(G)$ and $e(G)$ denote the number of vertices and the number of edges in $G$, respectively. By Euler's polyhedral formula, we have
\begin{equation}\label{3}
e(G)=v(G)+f(G)-2=v(G)+t-1.
\end{equation}
Let us call an edge a {\em boundary edge} if it belongs to the boundary of the polygon $P$. A boundary edge is said to be {\em full} if it is the full side of a triangle in $\mathcal T$, and {\em partial} if it is a proper part of a side. (Note that partial boundary edges did not show up in the proof of Theorem~6, because there $P$ was a {\em convex} polygon.) Let $e_{\rm full}$ and $e_{\rm part}$ stand for the number of full and partial boundary edges, respectively.

Note that all boundary edges are contained in the union of the boundaries of the triangles in ${\mathcal T}'$. The length of any side is at least $\epsilon'$. Thus, in view of (\ref{1}), the number of full boundary edges satisfies the inequality
\begin{equation}\label{4}
e_{\rm full}\le\frac{t'}{\epsilon'}=O(r).
\end{equation}
To bound the number of partial boundary edges, it is enough to observe that every partial boundary edge contains a vertex of a triangle $T\in{\mathcal T}'$, and any vertex of $T$ belongs to at most one partial boundary edge. Hence, we have
\begin{equation}\label{5}
e_{\rm part}\le3t'=O(r).
\end{equation}

A {\em stretch} is a {\em minimal} segment that can be decomposed into {\em sides} (of triangles in $\mathcal T$) and {\em partial boundary edges} in two different ways. (Note that in the proof of Theorem~6 we did not have to include boundary edges in the definition of a stretch, because in that setting all stretches lied in the interior of $P$.) The number of sides belonging to (contained in) a stretch is called the {\em size} of the stretch. Observe that sides contained in the boundary of $P$ (full boundary edges) do not belong to any stretch, but every other side belongs to a unique one. Thus, the total size of the stretches is equal $3t-e_{\rm full}$. We call a stretch {\em improper} if it contains a partial boundary edge. Otherwise, it is called a {\em proper} stretch. The size of a proper stretch is at least $3$, the size of an improper stretch is at least $2$. We call a proper stretch of size $3$ {\em tight}, and denote the number of tight stretches by $\sigma_{\rm tight}$. A stretch is said to be {\em loose} if it is not tight. Let $L_{\rm loose}$ denote the total length of all loose stretches.
\begin{figure}\centering
\begin{tikzpicture}[scale=0.6]


\draw[thick,dashed] (-3,-1) -- (6,2) -- (-3,8) -- (0,0);
\draw[thick,dashed] (-8,2) -- (-2,1) -- (-3,8);
\draw[thick] (-3,8) -- (-6,3);
\draw[thick,dashed] (-3,-1) -- (-1.05,2.9);
\draw[thick] (-8,2) -- (-2.5,4.75);
\draw[ultra thick,red] (-6,1.66) -- (-0.8,-3);
\draw[thick] (-0.8,-3)-- (-8,2);
\draw[ultra thick,red] (-1.9,-2) -- (10,4);
\draw[thick] (2,-2)  --(-1.9,-2);
\draw[thick] (-0.8,-3) -- (1,-2);
\draw[thick] (2,-0.05)--(2,-3)-- (10,4);
\draw[thick,dashed] (0,6)--(15,3);
\draw[thick] (15,3) -- (9,9)--(0,6);
\draw[thick] (6,0.5)--(15,3);
\draw[thick] (-3,8)--(6,8);
\node[fill=white, circle, inner sep = 0pt] at (-2.7,2.7) {s};
\node[fill=white, circle, inner sep = 0pt] at (-3,5.7) {s};
\node[fill=white, circle, inner sep = 0pt] at (-2.07,3.5) {l};
\node[fill=white, circle, inner sep = 0pt] at (-1.6,2.4) {s};
\node[fill=white, circle, inner sep = 0pt] at (-1.8,4.2) {s};
\node[fill=white, circle, inner sep = 0pt] at (-0.6,3) {l};
\node[fill=white, circle, inner sep = 0pt] at (-0.7,1) {s};
\node[fill=white, circle, inner sep = 0pt] at (-1.7,0.8) {l};
\node[fill=white, circle, inner sep = 0pt] at (-1.4,-0.1) {s};
\node[fill=white, circle, inner sep = 0pt] at (-4.7,2) {l};
\node[fill=white, circle, inner sep = 0pt] at (-3.8,1) {s};
\node[fill=white, circle, inner sep = 0pt] at (-6.6,1.4) {s};
\node[fill=white, circle, inner sep = 0pt] at (-2.9,0.1) {s};
\node[fill=white, circle, inner sep = 0pt] at (2,1) {s};
\node[fill=white, circle, inner sep = 0pt] at (0,-0.5) {l};
\node[fill=white, circle, inner sep = 0pt] at (1,4.8) {l};
\node[fill=white, circle, inner sep = 0pt] at (4,3.8) {s};
\node[fill=white, circle, inner sep = 0pt] at (-0.8,7) {s};
\node[fill=white, circle, inner sep = 0pt] at (8,5) {l};
\node[fill=white, circle, inner sep = 0pt] at (6,4.3) {s};
\node[fill=white, circle, inner sep = 0pt] at (12,3.1) {s};
\node[fill=white, circle, inner sep = 0pt] at (5.5,-1.5) {partial};
\node[fill=white, circle, inner sep = 0pt] at (11,1.1) {full};

\end{tikzpicture}
\caption{Part of a tiling without two triangles sharing an edge. The short and long parts of each tight stretch (dashed) are marked by s and l, respectively. Loose stretches  are marked thick and red. One partial and one full boundary edge is marked.}
\end{figure}
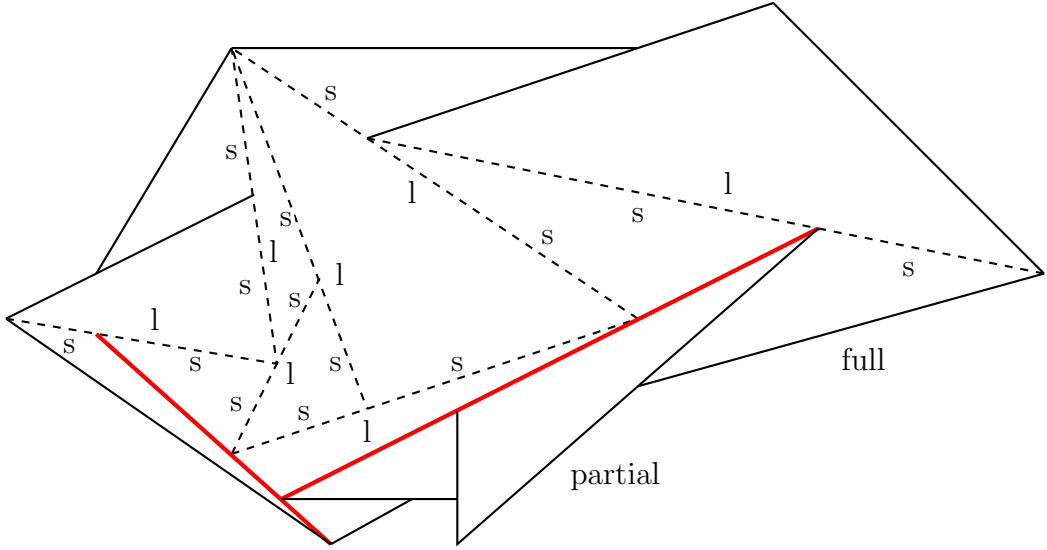

Summing up the lengths of tight and loose stretches, we obtain
$$3t-e_{\rm full}=3\sigma_{\rm tight}+L_{\rm loose},$$ whence
\begin{equation}\label{6}
t=\sigma_{\rm tight}+{L_{\rm loose}\over 3} + {e_{\rm full}\over 3}.
\end{equation}

A vertex is called {\em subdividing} if it is the interior point of a side. Let $v^*$ denote the number of subdividing vertices. A subdividing vertex is an interior point of a single side and, hence, of a single stretch. As in the proof of Theorem~6, the number of subdividing vertices in the interior of a proper stretch of size $s$ is exactly $s-2$. The number of subdividing vertices of an improper stretch of size $s$ is at least $s-1$. Thus, any loose stretch of size $s$ has at least $s/2$ subdividing vertices. Hence, the total number of subdividing vertices satisfies
\begin{equation}\label{7}
v^*\ge \sigma_{\rm tight}+{L_{\rm loose}\over 2}.
\end{equation}

The sum of the number of edges on all faces of $G$ is equal to $2e(G)$. The $t$ triangles have $3t+v^*$ edges in total, while the infinite face has $e_{\rm full}+e_{\rm part}$ edges (namely, the boundary edges). Hence,
\begin{equation}\label{8}
2e(G)=3t+v^*+e_{\rm full}+e_{\rm part}.
\end{equation}

Comparing this with (\ref{3}), we get $2v^*+2t-2\le 2e(G)=3t+v^*+e_{\rm full}+e_{\rm part},$ so that
$$v^*\le t+e_{\rm full}+e_{\rm part}+2.$$
If we plug into this inequality the estimates~(\ref{6}) and~(\ref{7}), we conclude that
\[L_{\rm loose}\le8e_{\rm full}+6e_{\rm part}+12.\]
Taking~(\ref4) and~(\ref5) into account, this implies
\begin{equation}\label{9}
L_{\rm loose}=O(r).
\end{equation}

Since the number of sides on loose stretches as well as the number of sides that are full boundary edges is $O(r)$, we should concentrate on the tight stretches containing most of the $3t=\Omega(r^2)$ sides. We call the longest side on any tight stretch {\em long} and the other two sides on the same tight stretch {\em short}. The sides not belonging to any tight stretch are neither long nor short.
\smallskip

Along each tight stretch, the length of the long side is equal to the sum of the lengths of the two short sides. Therefore, the total length of all short sides is the same as the total length of all long sides. Define the quantity $W$ as

\begin{equation*}
\begin{aligned}
W=&\left({\frac23}-\epsilon''\right)(\# \textrm{long sides})-{\frac13}(\# \textrm{short sides})\\
&+(\textrm{total length of short sides})-(\textrm{total length of long sides}).
\end{aligned}
\end{equation*}

Since each tight stretch has one long and two short sides, we have $\# \textrm{long sides}=\sigma_{\rm tight}$ and $\# \textrm{short sides}=2\sigma_{\rm tight}$. The last two terms of $W$ cancel out. Thus,
\begin{equation}\label{10}
W=-\epsilon''\sigma_{\rm tight}.
\end{equation}

On the other hand, we can compute $W$ by adding up the contributions of the triangles $T\in {\mathcal T}$. We classify the triangles in $\mathcal T$ as follows. If a triangle has a side that it neither short, nor long, we call it {\em exceptional}. If $T$ is not exceptional and it has $i$ long sides, we say that $T$ is of {\em type} $i$, for $i=0,1,2,3$. If $T$ is of
\smallskip

\noindent{\bf type 0}, its contribution to $W$ is $0$, as the total length of its short sides is the perimeter of $T$, which is equal to $\frac13(\# \textrm{short sides of T})=1$.

\noindent{\bf type 1}, its contribution to $W$ is the total length of its two short sides minus the length of its long side minus $\epsilon''$, which is nonnegative, according to the definition of $\epsilon''$.

\noindent{\bf type 2}, its contribution to $W$ is $2(\textrm{length of its short side})-2\epsilon''\ge2\epsilon'-2\epsilon''\ge0$.

\noindent{\bf type 3}, its contribution to $W$ is equal to $1-3\epsilon''>0$.
\smallskip

\smallskip
The only triangles whose contribution to $W$ can be negative are the exceptional ones. The contribution of each exceptional triangle is is at least $-\frac23$. The number of these triangles in $\mathcal T$ is at most $L_{\rm loose}+e_{\rm full}$, so by (\ref4) and (\ref9) we obtain $-W=O(r)$.
Comparing this bound to (\ref{10}), we have
\[\sigma_{\rm tight}=O(r).\]
Plugging (\ref{4}), (\ref{9}), and this bound into~(\ref6), we conclude that $t=O(r)$, which contradicts~(\ref1). The contradiction proves the theorem. $\Box$

\section{Triangles with similar sides---Proof of Theorem 5}

Theorem 5 is an easy corollary to the following result.
\medskip

\noindent{\bf Theorem 7.} {\em In any locally finite tiling of the plane with triangles, there is a triangle with a side that is the union of one or more sides of other triangles.}
\medskip

\noindent{\bf Proof.} Assume we have a tiling $\mathcal T$ contradicting the statement of the theorem, and let $pst$ be any triangle in $\mathcal T$. As $pt$ cannot be obtained as the union of sides of other triangles, $p$ or $t$ must be an interior point of a side of another triangle along the line $pt$. Similar statements are true for the other two sides of the triangle $pst$, but no vertex can be an interior point of more than one side in the tiling. Thus, there are three triangles along the lines $pt$, $ts$, and $sp$, that contain $p$, $t$, and $s$, respectively, in this {\em cyclic} order. Similar statements hold for all triangles in the tiling.
\smallskip

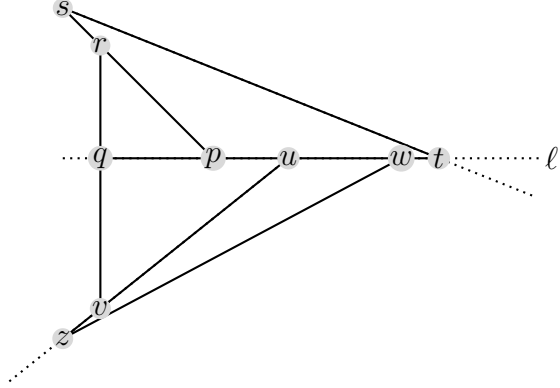
\begin{figure}\centering
\begin{tikzpicture}[scale=1]


\draw[thick,dotted] (-4,0) -- (2.5,0);
\draw[thick] (-2,0) -- (1,0)-- (-4,2)--cycle;
\draw[thick] (-3.5,1.5) -- (-3.5,0)-- (-2,0);
\draw[thick,dotted] (-4,2) -- (2.25,-0.5);
\draw[thick] (-3.5,0) -- (-3.5,-2)-- (-1,0);
\draw[thick,dotted] (-1,0) -- (-4.75,-3);
\draw[thick] (-3.5,-2) -- (-4,-2.4)-- (0.5,0);

\node[fill=white, circle, inner sep = -0pt, minimum size=0pt] at (2.5,0) {$\ell$};

\node[fill=Grey!30, circle, inner sep = -0pt, minimum size=0pt] at (-2,0) {$p$};
\node[fill=Grey!30, circle, inner sep = -0pt, minimum size=0pt] at (-4,2) {$s$};
\node[fill=Grey!30, circle, inner sep = -0pt, minimum size=0pt] at (1,0) {$t$};
\node[fill=Grey!30, circle, inner sep = -0pt, minimum size=0pt] at (-3.5,1.5) {$r$};
\node[fill=Grey!30, circle, inner sep = -0pt, minimum size=0pt] at (-3.5,0) {$q$};
\node[fill=Grey!30, circle, inner sep = -0pt, minimum size=0pt] at (-3.5,-2) {$v$};
\node[fill=Grey!30, circle, inner sep = -0pt, minimum size=0pt] at (-1,0) {$u$};
\node[fill=Grey!30, circle, inner sep = -0pt, minimum size=0pt] at (0.5,0) {$w$};
\node[fill=Grey!30, circle, inner sep = -0pt, minimum size=0pt] at (-4,-2.4) {$z$};


\end{tikzpicture}
\caption{Illustration for the proof of Theorem 7.}\label{fig5v2}
\end{figure}

We have seen that there is a unique triangle $T\in \mathcal T$ that has $p$ as an interior point of one of its sides. Assume without loss of generality that this side of $T$ belongs to the line $\ell$ containing $pt$. Notice that there is precisely one other triangle in $\mathcal T$, different from $pst$, whose vertex is $p$, and this triangle must also have a side contained in $\ell$. (Indeed, if there were a triangle in $\mathcal T$ with $p$ as a vertex that does not have a side contained in $\ell$, it would not satisfy the cyclic property described in the previous paragraph.) In the same way, for any other vertex $p'$ of the tiling, there are precisely two triangles with this vertex, and there is a line that contains one side of each. Let $pqr$ and $pst$ be the two triangles containing $p$ as a vertex, where $p$, $r$, and $s$ are collinear; see Figure~\ref{fig5v2}.







Since $p$ is an interior point of a side along $\ell$, we conclude that $q$ must be an interior point of a side along the line $qr$ and $t$ must be an interior point of a side along the line $st$. Let $quv$ be the other triangle with vertex $q$ along $\ell$, where $v$ lies on the line $qr$. As $q$ is not an interior point of a side along the line $\ell$, so $u$ must be an interior point of a side along $\ell$: either of the side $qp$ or of $pt$ (as the line $\ell$ is ``blocked'' by another triangle at $t$). Finally, consider the other triangle of the tiling that has $u$ as a vertex. It must have another vertex $w$ along $\ell$ (in fact, in the segment $ut$, because $\ell$ is ``blocked'' at $t$) and a third vertex $z$ on the ray emanating from $u$ toward $v$ (see Figure~\ref{fig5v2}).
Now $u$ is the interior point of a side along $\ell$, so $w$ must be an interior point of a side along the line $wz$. This is possible only if $w=t$ and $z$ lies on the line $st$, because this is the first line ``blocking'' $\ell$. The ray from $u$ toward $v$ (containing $z$) does not intersect the line $st$ (which is also supposed to contain $z$). This contradiction proves the theorem.
$\Box$

\medskip

Two triangles of a locally finite triangular tiling $\mathcal T$ are called {\em neighbors} if they share a boundary segment. It follows from the proof of the last theorem that for any triangle $T\in \mathcal T$, either $T$ itself, or one of its neighbors, second neighbors, or third neighbors has a side that is the union of one or more sides of other triangles.

\end{document}